\input amstex
\documentstyle{amsppt}
\topmatter
\NoRunningHeads
\title
On the paper
"The Payne conjecture for Dirichlet and buckling eigenvalues", by Genqian Liu
\endtitle
\author
Leonid Friedlander
\endauthor
\affil
University of Arizona
\endaffil
\email
friedlan@math.arizona.edu
\endemail

\endtopmatter
\document
Let $\Omega$ be a bounded domain in $\Bbb R^n$ with smooth boundary. Let $\lambda_1<\lambda_2\leq\cdots$ be the eigenvalues of the Dirichlet Laplacian in $\Omega$, and let $\Lambda_1\leq\Lambda_2\leq\cdots$ be the eigenvalues of the buckling problem: these are the values of $\Lambda$ for which the equation
$$\Delta^2u+\Lambda\Delta u=0\quad\text{in}\ \Omega,\tag 1$$
with the Diriclet boundary conditions $u=(\partial u/\partial n)=0$ on $\partial\Omega$ has a non-trivial solution.  Here $n$ is the field of unit outward normal vectors to $\partial\Omega$. The conjecture that G. Liu is trying to settle in [L] is that  $\Lambda_k\geq \lambda_{k+1}$ for every $k$. Unfortunately, the proof contains a mistake that looks irreparable.

Let me very briefly explain the nature of the problem and what is wrong in [L]. This conjectue is somewhat similar to the inequality that I proved in [F]. Let $0=\mu_1<\mu_2<\cdots$ be the eigenvalues of the Neumann Laplacian in $\Omega$. I proved that $\lambda_k\geq \mu_{k+1}$. The proof used the Dirichlet-to-Neumann operator $R(\lambda)$. It is defined when $\lambda$ is different from $\lambda_k$'s, and it takes a function $\phi(x)$ that is defined on $\partial\Omega$ to $\partial u/\partial n$ where $u(x)$ is the solution  to the problem
$$\Delta u+\lambda u=0\ \text{in}\ \Omega,\quad u(x)=\phi(x)\ \text{on}\ \partial\Omega.$$
It turns out that $R(\Lambda)$ is a semibounded operator with discrete spectrum, and the number of its negative eigenvalues equals $N_N(\lambda)-N_D(\lambda)$ whare $N_*(\lambda)$ is the number of eigenvalues of the corresponding problem  that are smaller than $\lambda$. With that, proving $\lambda_k\geq\mu_{k+1}$ becomes equivalent to showing that the operator $R(\lambda)$ has at least one negative eigenvalue when $\lambda>0$. G. Liu is trying to go along similar lines.  The eigenvalue problem for the Dirichlet Laplacian is equivalent to (1), with the boundary conditions $u(x)=\Delta u(x)=0$ on $\partial\Omega$. The Dirichlet-to-Neumann operator can be replaced by the Neumann-to-Laplacian operator $S(\Lambda)$ that sends a function
$\psi(x)$  on $\partial\Omega$ to the restriction of $\Delta u$ to $\partial\Omega$ where $u(x)$ is the solution of (1) in $\Omega$, wth the boundary conditions $u=0$, $(\partial u/\partial n)=\psi$ on $\partial\Omega$. The operator $S(\Lambda)$ is defined when $\Lambda\not=\Lambda_j$.  Liu proves that the number of negative eigenvalues of $S(\lambda)$ equals $N_D(\lambda)-N_{\text{buck}}(\lambda)$. Then, to prove Payne's conjecture one has to show that $S(\lambda)$ has a negative eigenvalue for $\lambda>\lambda_1$.

Let $\beta_1(\lambda)$ be the smallest eigenvalue of $S(\lambda)$. Then $\beta_1(\lambda)$ is the infimum of the Rayleigh quotient $(S(\lambda)\psi,\psi)/(\psi,\psi)$. This Rayleigh quotient can be re-written in terms of the function $v(x)$  that solves the problem $\Delta^2 v+\lambda\Delta v=0$ in $\Omega$, $v=0$, $(\partial v/\partial n)=\psi$ on $\partial\Omega$:
$$\beta_1(\lambda)=\inf_{v\in\Cal L(\lambda)}\frac{\int_\Omega(|\Delta v|^2-\lambda|\nabla v|^2)dx}{\int_{\partial\Omega}(\partial v/\partial n)^2dS}. \tag 2$$
Here $\Cal L(\lambda)$ is the space of functions such  that $\Delta^2 v+\lambda\Delta v=0$ in $\Omega$ and $v=0$ on $\partial\Omega$. The mistake in G. Liu's paper is that he replaced the space $\Cal L(\lambda)$ in (2) by a bigger space $H^2(\Omega)\cap H_0^1(\Omega)$ (formula (4.4) in [L].) It can actually be dome if $\lambda<\Lambda_1$. Then the quotient in (2) is bounded from below in $H^2(\Omega)\cap H_0^1(\Omega)$, the infimum is attained, and a minimizer belongs to $\Cal L(\lambda)$. The situation is drastically different for $\lambda>\Lambda_1$. In this case, this quotient is not bounded from below in $H^2(\Omega)\cap H_0^1(\Omega)$. To see that, take the ground state $u_1(x)$ of the buckling problem. One has
$$\alpha=\int_\Omega (|\Delta u(x)|^2-\lambda |\nabla(x)|^2)dx =-(\lambda-\Lambda_1)\int_\Omega |\nabla u(x)|^2dx<0.$$
Let $h(x)$ be any smooth function in $\bar\Omega$  that vanishes on $\partial\Omega$, and such that $\partial h/\partial n$ does not vanisih on $\partial\Omega$. Let $v_\epsilon=u+\epsilon h$.
Then the quotient in (2) computed for $v_\epsilon$ goes to $-\infty$ as $\epsilon\to 0$: the numerator approaches a negative number $\alpha$, the denominator goes to $0$.

Let me make three additional remarks.\newline
1. The main difficulty proving that $\beta_1(\lambda)<0$ is that test functions in (2) must satisfy the differential equation AND vanish on the boundary of $\Omega$. In the case of the Dirichlet-to-Neumann operator, there is no boundary condition for  test functions.\newline
2. The operator $S(\lambda)$ is positive when $\lambda<\lambda_1$, so a proof of $\beta_1(\lambda)<0$ should use the condition $\lambda>\lambda_1$ in a substantial way.\newline
3. The fact that the number of negative eigenvalues of the Dirichlet-to-Neumann operator (Neumann-to-Laplacian operator) equals the difference between counting functions of the spectra of corresponding operators is very general. It holds for any compact Riemannian manifold with boundary. Rafe Mazzeo proved in [M]  the inequality $\lambda_k\geq \mu_{k+1}$ for bounded domains in Riemannian symmetric spaces of noncompact type. He also indicated that positive curvature could pose a problem.  It  definitely fails for $\Bbb S^n_\epsilon$, which is
the $n$-dimensional sphere, with a cap of radius $\epsilon$ removed, if $\epsilon$ is small enough. I suspect, though I do not have a good  intuition here, that similar thing happens to the buckling problem. That would be interesting to know the asymptotics of $\Lambda_1(\Cal S^n_\epsilon)$ as $\epsilon\to 0$. If this is actually the case that the inequality $\beta_1(\lambda)<0$ when $\lambda>\lambda_1$ does not hold for any compact Riemannian manifold with bounday then its proof should use the specifics of an Euclidean space.

I wrote to Genquian Liu on September 7, 2021, communicated to him what  the mistake is, with rather detailed explanations. He did receive my letters. Nevertheless, his paper has not been withdrawn from the arxiv yet. So I decided to write this note.

\enddocument